\documentclass[reqno,12pt,a4paper]{amsart}

\usepackage{amssymb}
\usepackage{amsfonts}
\usepackage{latexsym}
\usepackage{amsthm}
\usepackage{bm}

\begin{document}

%%%%%%%%%% Some definitions %%%%%%%%%%

%%%%% Equations, labels, theorems %%%%%
\renewcommand{\theequation}{\arabic{section}.\arabic{equation}}
\def\lb{\label}                       \newcommand{\er}[1]{{\rm(\ref{#1})}}
\theoremstyle{plain}
\newtheorem{theorem}{\bf Theorem}[section]
\newtheorem{lemma}[theorem]{\bf Lemma}
\newtheorem{corollary}[theorem]{\bf Corollary}
\newtheorem{proposition}[theorem]{\bf Proposition}
\newtheorem{remark}[theorem]{\it Remark}
%\theoremstyle{remark}
%\newtheorem{remark}[theorem]{\bf Remark}

%%%%% Alphabet %%%%%
\def\a{\alpha}  \def\cA{{\mathcal A}}     \def\bA{{\bf A}}  \def\mA{{\mathscr A}}
\def\b{\beta}   \def\cB{{\mathcal B}}     \def\bB{{\bf B}}  \def\mB{{\mathscr B}}
\def\g{\gamma}  \def\cC{{\mathcal C}}     \def\bC{{\bf C}}  \def\mC{{\mathscr C}}
\def\G{\Gamma}  \def\cD{{\mathcal D}}     \def\bD{{\bf D}}  \def\mD{{\mathscr D}}
\def\d{\delta}  \def\cE{{\mathcal E}}     \def\bE{{\bf E}}  \def\mE{{\mathscr E}}
\def\D{\Delta}  \def\cF{{\mathcal F}}     \def\bF{{\bf F}}  \def\mF{{\mathscr F}}
\def\c{\chi}    \def\cG{{\mathcal G}}     \def\bG{{\bf G}}  \def\mG{{\mathscr G}}
\def\z{\zeta}   \def\cH{{\mathcal H}}     \def\bH{{\bf H}}  \def\mH{{\mathscr H}}
\def\e{\eta}    \def\cI{{\mathcal I}}     \def\bI{{\bf I}}  \def\mI{{\mathscr I}}
\def\p{\psi}    \def\cJ{{\mathcal J}}     \def\bJ{{\bf J}}  \def\mJ{{\mathscr J}}
\def\vT{\Theta} \def\cK{{\mathcal K}}     \def\bK{{\bf K}}  \def\mK{{\mathscr K}}
\def\k{\kappa}  \def\cL{{\mathcal L}}     \def\bL{{\bf L}}  \def\mL{{\mathscr L}}
\def\l{\lambda} \def\cM{{\mathcal M}}     \def\bM{{\bf M}}  \def\mM{{\mathscr M}}
\def\L{\Lambda} \def\cN{{\mathcal N}}     \def\bN{{\bf N}}  \def\mN{{\mathscr N}}
\def\m{\mu}     \def\cO{{\mathcal O}}     \def\bO{{\bf O}}  \def\mO{{\mathscr O}}
\def\n{\nu}     \def\cP{{\mathcal P}}     \def\bP{{\bf P}}  \def\mP{{\mathscr P}}
\def\r{\rho}    \def\cQ{{\mathcal Q}}     \def\bQ{{\bf Q}}  \def\mQ{{\mathscr Q}}
\def\s{\sigma}  \def\cR{{\mathcal R}}     \def\bR{{\bf R}}  \def\mR{{\mathscr R}}
\def\S{\Sigma}  \def\cS{{\mathcal S}}     \def\bS{{\bf S}}  \def\mS{{\mathscr S}}
\def\t{\tau}    \def\cT{{\mathcal T}}     \def\bT{{\bf T}}  \def\mT{{\mathscr T}}
\def\f{\phi}    \def\cU{{\mathcal U}}     \def\bU{{\bf U}}  \def\mU{{\mathscr U}}
\def\F{\Phi}    \def\cV{{\mathcal V}}     \def\bV{{\bf V}}  \def\mV{{\mathscr V}}
\def\P{\Psi}    \def\cW{{\mathcal W}}     \def\bW{{\bf W}}  \def\mW{{\mathscr W}}
\def\o{\omega}  \def\cX{{\mathcal X}}     \def\bX{{\bf X}}  \def\mX{{\mathscr X}}
\def\x{\xi}     \def\cY{{\mathcal Y}}     \def\bY{{\bf Y}}  \def\mY{{\mathscr Y}}
\def\X{\Xi}     \def\cZ{{\mathcal Z}}     \def\bZ{{\bf Z}}  \def\mZ{{\mathscr Z}}
\def\O{\Omega}

\def\ve{\varepsilon}   \def\vt{\vartheta}    \def\vp{\varphi}    \def\vk{\varkappa}

\def\Z{{\mathbb Z}}    \def\R{{\mathbb R}}   \def\C{{\mathbb C}}
\def\T{{\mathbb T}}    \def\N{{\mathbb N}}   \def\dD{{\mathbb D}}

%%%%% Arrows %%%%%

\def\la{\leftarrow}              \def\ra{\rightarrow}            \def\Ra{\Rightarrow}
\def\ua{\uparrow}                \def\da{\downarrow}
\def\lra{\leftrightarrow}        \def\Lra{\Leftrightarrow}

%%%%% Math signs %%%%%

\let\ge\geqslant                 \let\le\leqslant
\def\/{\over}                    \def\iy{\infty}
\def\sm{\setminus}               \def\es{\emptyset}
\def\ss{\subset}                 \def\ts{\times}
\def\pa{\partial}                \def\os{\oplus}
\def\ev{\equiv}                  \def\1{1\!\!1}
\def\iint{\int\!\!\!\int}        \def\iintt{\mathop{\int\!\!\int\!\!\dots\!\!\int}\limits}
\def\el2{\ell^{\,2}}

%%%%% Typography %%%%%

\def\lt{\biggl}                  \def\rt{\biggr}
\def\no{\noindent}               \def\ol{\overline}
\def\BBox{\hspace{1mm}\vrule height6pt width5.5pt depth0pt \hspace{6pt}}
\def\wt{\widetilde}              \def\wh{\widehat}
\newcommand{\nt}[1]{{\mathop{#1}\limits^{{}_{\,\bf{\sim}}}}\vphantom{#1}}
\newcommand{\nh}[1]{{\mathop{#1}\limits^{{}_{\,\bf{\wedge}}}}\vphantom{#1}}
\newcommand{\nc}[1]{{\mathop{#1}\limits^{{}_{\,\bf{\vee}}}}\vphantom{#1}}
\newcommand{\oo}[1]{{\mathop{#1}\limits^{\,\circ}}\vphantom{#1}}
\newcommand{\po}[1]{{\mathop{#1}\limits^{\phantom{\circ}}}\vphantom{#1}}

%%%%% Math operations %%%%%

\def\Im{\mathop{\rm Im}\nolimits}
\def\Iso{\mathop{\rm Iso}\nolimits}
\def\Ker{\mathop{\rm Ker}\nolimits}
\def\Ran{\mathop{\rm Ran}\nolimits}
\def\Re{\mathop{\rm Re}\nolimits}
\def\Tr{\mathop{\rm Tr}\nolimits}
\def\arg{\mathop{\rm arg}\nolimits}
\def\const{\mathop{\rm const}\nolimits}
\def\det{\mathop{\rm det}\nolimits}
\def\diag{\mathop{\rm diag}\nolimits}
\def\dim{\mathop{\rm dim}\nolimits}
\def\dist{\mathop{\rm dist}\nolimits}
\def\res{\mathop{\rm res}\limits}
\def\sign{\mathop{\rm sign}\nolimits}
\def\supp{\mathop{\rm supp}\nolimits}

\def\cF{{\mathcal F}}

%%%%%%%%%%% End of definitions %%%%%%%%%%

\title[Application of fixed point theorem to inverse problem]
{An application of the fixed point theorem\\ to the inverse Sturm-Liouville problem}

\author{Dmitry Chelkak}

\dedicatory{To Nina N. Uraltseva on the occasion of her 75th birthday}

%\subjclass{34B24, 34A55}

%\keywords{Sturm-Liouville operators, characterization of spectral data}

\thanks{\textsc{Dept. of Math. Analysis, St.~Petersburg State University.
28~Universitetskij~pr., Staryj Petergof, 198504 St.~Petersburg, Russia.} \quad Partly funded
by the RF President grants MK-4306.2008.1, NSh-2409.2008.1, and P.Deligne's 2004 Balzan prize
in Mathematics.}

\begin{abstract}
We consider Sturm-Liouville operators $-y''+v(x)y$ on $[0,1]$ with Dirichlet boundary
conditions $y(0)=y(1)=0$. For any $1\le p<\infty$, we give a short proof of the
characterization theorem for the spectral data corresponding to $v\in L^p(0,1)$.
\end{abstract}

\maketitle

\section{Introduction}

In this paper we consider the inverse spectral problem for self-adjoint Sturm-Liouville
operators
\begin{equation}
\label{SLOp} \cL y=-y''+v(x)y,\qquad y(0)=y(1)=0,
\end{equation}
acting in the Hilbert space $L^2(0,1)$, with $v\in L^p(0,1)$ for some (fixed) $1\!\le\!
p\!<\!\infty$ (we denote by $\|v\|_p<\infty$ the standard $L^p$ norm of $v$). The spectrum of
$\cL$ is denoted by
\[
\l_1(v)<\l_2(v)<\l_3(v)<\dots
\]
It is purely discrete, simple and satisfies the asymptotics
\[
\l_n(v)=\pi^2 n^2 + \nh v^{(0)} + \mu_n(v),
\]
where
\[
\textstyle \nh v^{(0)}:= \int_0^1v(t)dt\qquad \mathrm{and} \qquad \mu_n(v)=o(1)\ \
\mathrm{as}\ \ n\to\infty.
\]
Note that $\nh v^{(0)}$ can be immediately reconstructed from the Dirichlet spectrum as the
leading term in the asymptotics of $\l_n(v)-\pi^2n^2$.

Starting with the famous uniqueness theorem of Borg \cite{Bo}, the inverse spectral theory of
scalar 1D differential operators was developed in detail, and currently there are several
classical monographs devoted to different approaches to these problems (see,
e.g.~\cite{MaBook}, \cite{LeBook}, \cite{PT}). Traditionally, the principal attention is paid
to explicit reconstruction procedures that allow one to find the unknown potential starting
from a given spectral data. The careful analysis of these procedures (a)~has the significant
practical interest; \mbox{(b)~allows} one to find the necessary and sufficient conditions for
spectral data to correspond to some potential from a given class. The latter results are
usually called {\it characterization theorems}. In order words, they say that the mapping
\[
\cM\ :\ \{\mathrm{potentials}~v(x)\}\to\{\mathrm{spectral\ data}\}
\]
is a bijection between some fixed space of potentials $\cE$ and a class of spectral data $\cS$
which is described explicitly.

Our main result is a {\it short proof} of characterization theorems for spectral data of
Sturm-Liouville operators (\ref{SLOp}) corresponding to $L^p$ potentials (see
Theorem~\ref{CharEven} and Theorems~\ref{CharMarchenko},~\ref{CharTrubowitz} below). At least
for $p=1$ and $p=2$ these results are well known in the literature, but, unfortunately, we
don't know any reference that cover all $p$'s simultaneously. It is worthwhile to emphasize
that the main goal of our paper is to present the {\bf method} (more precisely, the {\it
simplification} of Trubowitz's scheme, see below) rather than new results. We hope that this
method is applicable to other inverse spectral problems too.

\smallskip

For simplicity, we first focus on symmetric (or {\it even}) potentials
\[
v(x)\equiv v(1-x),\qquad x\in [0,1].
\]
Then, it is well known that the spectrum itself determines a potential uniquely (see,
e.g.~\cite{PT}, pp. 55--57 and p. 62 for a very short proof). Let
\begin{equation}
\label{MapEven} \cM\ :\ v\mapsto (\,\nh v^{(0)}\,,\, \{\,\mu_n(v)\,\}_{n=1}^{\infty}\,)
\end{equation}
and
\[
\textstyle \cF_\mathrm{cos}: v\mapsto \left(\nh v^{(0)}; \{-\nh
v^{(cn)}\}_{n=1}^{\infty}\right),\quad\mathrm{where}\quad \nh v^{(cn)}= \int_0^1 v(t)\cos
(2\pi n t) dt,
\]
denote (up to a sign) the cosine-Fourier transform.

\begin{theorem}
\label{CharEven} Let $1\le p<\infty$. The mapping $\cM$ given by (\ref{MapEven}) is a
bijection between the space of all symmetric $L^p$-potentials $\cE=L^p_\mathrm{even}(0,1)$ and
the subset $\cS$ of the Fourier image $\cF_\mathrm{cos}L^p_\mathrm{even}(0,1)$ consisting of
all sequences $ \mu^*=(\mu^*_0,\mu^*_1,\mu^*_2,\dots)\in
\cF_\mathrm{cos}L^p_\mathrm{even}(0,1) $ such that
\begin{equation}
\label{MuLe}
\pi^2n^2+\mu^*_0+\mu^*_n<\pi^2(n\!+\!1)^2+\mu^*_0+\mu^*_{n+1}\quad\mathit{for~all}~n\ge 1.
\end{equation}
\end{theorem}

In general, in order to prove the characterization theorem, one needs

\smallskip

\noindent (i)\phantom{ii} to solve the direct problem, i.e., to show that $\cM$ maps $\cE$
into $\cS$;

\noindent (ii)\phantom{i} to prove the uniqueness theorem, i.e., the fact that the mapping
$\cM$ is $1$-to-$1$;

\noindent (iii) to prove that $\cM$ is a surjection.

\smallskip

Usually, the first part is rather straightforward and the second part can be simply done in a
nonconstructive way without any references to explicit reconstruction procedures. Thus, the
hardest part of such theorems is the third one. It was suggested by Trubowitz and co-authors
(see~\cite{PT}) to use the following abstract scheme in order to prove (iii). To omit
inessential technical details concerning the particular structure of the infinite-dimensional
manifold $\cS$ (the restriction (\ref{MuLe}) in our case), in the next paragraph we think of
$\cS$ as of a Banach space equipped with the usual addition operation. Following Trubowitz's
scheme, it is sufficient

\smallskip

\noindent (a) to show that $\cM(\cE)$ contains some open set $\cO\ss\cS$ (say, some
neighborhood of $0$);

\noindent (b) to show that for some dense subset $\cL\ss\cS$ the following is fulfilled:

\begin{center}for any $s\in\cM(\cE)$ and $l\in\cL$, one has $s+l\in\cM(\cE)$.\end{center}

\noindent Since, for any $s\in\cS$, the set $s-\cL$ is dense in $\cS$, one has $s=o+l$ for
some $o\in\cO$ and $l\in\cL$. Thus, (b) implies $s\in\cM(\cE)$ because, due to (a),
$o\in\cO\ss\cM(\cE)$.

\smallskip

Loosely speaking, to prove (b) one needs to apply the reconstruction procedure only for {\it
``nice'' perturbations} $l\in\cL$ of spectral data (but starting with an arbitrary $s\in\cO$).
Following \cite{PT}, (a) can be deduced from the implicit function theorem applied to the
mapping $\cM$ near $v=0$. In order to do this, it is necessary to prove that $\cM$ is
continuously differentiable (in appropriate spaces) {\it everywhere near $v=0$}. Actually, the
proving of the differentiability of $\cM$ near $0$ is not much simpler than the
differentiability of $\cM$ {\it everywhere in $\cE$}, since the information about the norm
$\|v\|$ doesn't help to prove the existence of the Fr\'echet derivative $d_v\cM$ at $v$. The
main purpose of our paper is to point out that, in fact, one can notably simplify this part of
the proof, using some abstract fixed point theorem and

\smallskip

\noindent (a1)~the differentiability of $\cM$ (in the Fr\'echet sense) at {\it only one point
$v=0$};

\noindent (a2) the continuity of $\cM$ in the weak-$*$ topology (if $1<p<\infty$).

\smallskip

The paper is organized as follows. We start with some preliminaries in Sect.~\ref{SectPre}.
The very simple but crucial application of the Leray-Schauder-Tyhonoff fixed point theorem
which allows us to (almost immediately) derive (a) from (a1) and (a2) is given in
Sect.~\ref{SectLocSurjCore}. The properties (a1), (a2), and (a) for the mapping
(\ref{MapEven}) are proved in Sect.~\ref{SectLocalSurjLp}, if $1<p<\infty$. The necessary
modifications for $p=1$ are given in Sect.~\ref{SectL1}. The proof of Theorem~\ref{CharEven}
is finished in Sect.~\ref{SectGlob}. For the sake of completeness, in Sect.~\ref{SectGen} we
also consider nonsymmetric potentials. For both usual choices of additional spectral data
(Marchenko's normalizing constants as well as Trubowitz'z norming constants), we prove the
characterization Theorems~\ref{CharTrubowitz}~and~\ref{CharMarchenko} similar to
Theorem~\ref{CharEven}.

\smallskip

The scheme described above is quite general and can be used to prove similar characterization
theorems for other ``reasonable'' spaces of potentials instead~of~$L^p(0,1)$. Another approach
to these results (for $W^\theta_2$ potentials with $\theta\ge -1$) based on the interpolation
technique was suggested in \cite{SS08}.

\smallskip

\noindent {\bf Acknowledgements.} It is my pleasure to dedicate this paper to Nina
\mbox{N.~Uraltseva}, whose lectures on PDE I had a chance to attend, as many other generations
of students. The author is also grateful to Evgeny Korotyaev, Boris M.~Makarov and Sasha
Pushnitski for helpful discussions.

\section{Symmetric case, proof of Theorem \ref{CharEven}}
\setcounter{equation}{0}

\subsection{Preliminaries}
\label{SectPre}

Let $\vp(x,\l,v)$ denote the solution to the differential equation
\begin{equation}
\label{DiffEq} -y''+v(x)y=\l y
\end{equation}
satisfying the initial conditions $\vp(0,\l,v)=0$, $\vp'(0,\l,v)=1$. It can be constructed by
iterations as
\begin{equation}
\label{Vp=Series} \vp(x,\l,v)=\sum_{k=0}^{\infty}\vp_k(x,\l,v),\qquad \mathrm{where}\qquad
\vp_0(x,\l)={\sin\sqrt{\l}x}\big/{\sqrt{\l}}
\end{equation}
and
\begin{align}
\label{VpK=}
\vp_k(x,\l,q) & =\int_0^x\vp_0(x\!-\!t,\l)\vp_{k-1}(t,\l,v)v(t)dt \\
 &=\iintt_{0=t_0\le t_1\le\dots \le t_k\le t_{k+1}=x} {\textstyle \prod_{m=0}^k
\vp_0(t_{m+1}\!-\!t_m,\l)}\cdot v(t_1)\dots v(t_k) dt_1\dots dt_k\,. \notag
\end{align}
Since $|\vp_0(t,\l)|\le e^{|\Im\sqrt{\l}|t}/\sqrt{\lambda}$, one immediately obtains the
estimate
\begin{equation}
\label{VpkEstimate} |\vp_k(x,\l,v)| \le \frac{\|v\|_1^k} {k!}\cdot
\frac{e^{|\Im\sqrt{\l}|x}}{|\l|^{(k+1)/2}}\,.
\end{equation}
In particular, the series (\ref{Vp=Series}) converges uniformly in $\l$ and $v$ on bounded
subsets of $\C$ and $L^1(0,1)$, respectively. Since $\l_n(v)$ are the zeros of the entire
function
\[
w(\l,v):=\vp(1,\l,v)=\sum_{k=0}^{\infty}\vp_k(1,\l,v),\qquad \l\in\C,
\]
and the zeros of $\vp_0(1,\l)$ are $\pi^2n^2$, (\ref{VpkEstimate}) easily gives
$\l_n(v)=\pi^2n^2+O(\|v\|_1 e^{\|v\|_1})$. Taking into account the second term
$\vp_1(1,\l,v)$, one obtains
\begin{equation}
\label{Lasympt} \l_n(v)=\pi^2n^2+ \nh v^{(0)} - \nh v^{(cn)} + O\lt(\frac{\|v\|^2_1\,
e^{\|v\|_1}}{n}\rt),
\end{equation}
with some absolute constant in the $O$-bound. We also need the following simple Lemma
\begin{lemma}
\label{WeakCont} Let $1<p<\infty$ and $v_s,v\in L^p(0,1)$ be such that $v_s\to v$ weakly in
$L^p(0,1)$ as $s\to\infty$. Then $\l_n(v_s)\to\l_n(v)$ for any $n\ge 1$.
\end{lemma}

\begin{proof} Cf.~\cite{PT}~p.~18.
Since $\l_n(s)$ are the zeros of the entire functions $w(\cdot,v_s)$, it is sufficient to
prove that $w(\l,v_s)\to w(\l,v)$ uniformly in $\l$ on bounded subsets of $\C$. Let
$q=p/(p-1)$. For any $k\ge 1$, the functions
\[
f_{\l}(t_1,t_2,...,t_k):=\textstyle \chi_{\{0=t_0\le t_1\le \dots \le t_k\le
t_{k+1}=1\}}\cdot\prod_{m=0}^k\vp_0(t_{m+1}-t_m,\l),\qquad |\l|\le M,
\]
form a compact set in $L^q([0,1]^k)$. Thus, since $\prod_{m=1}^k v_s(t_m)\to \prod_{m=1}^k
v_s(t_m)$ weakly in $L^p([0,1]^k)$, one has $\vp_k(1,\l,v_s)\to \vp_k(1,\l,v)$ uniformly in
$\l:|\l|\le M$. As the norms $\|v_s\|_p$ are uniformly bounded and the series
(\ref{Vp=Series}) converges uniformly in $\l$ and $v$ on bounded subsets, it implies
$w(\l,v_s)\to w(\l,v)$ uniformly in $\l$ on bounded subsets.
\end{proof}

\subsection{Local surjection near $\bm{v=0}$. Core argument.}
\label{SectLocSurjCore}

\begin{lemma}
\label{LocSurjLemma} Let $E$ be a reflexive Banach space and a mapping $\Phi:B_E(0,r)\to E$ be
defined in some neighborhood $B_E(0,r)=\{v\in E:\|v\|_E<r\}$ of $0$. If $\Phi$ is

\smallskip

(a1) differentiable in the Fr\'echet sense at $0$ and $d_0\Phi=I$, i.e.,
\[
\|\Phi(v)-v\|_E=o(\|v\|_E)~\mathit{as}~\|v\|_E\to 0;
\]

(a2) continuous in the weak topology, i.e.,
\[
v_n\to v~\mathit{weakly}~~~\Rightarrow~~~\Phi(v_n)\to\Phi(v)~\mathit{weakly},
\]

\noindent then $\Phi$ is a local surjection at $0$, i.e., $\Phi(B_E(0,r))\supset
B_E(0,\delta)$ for some $\delta>0$.
\end{lemma}

\begin{proof}
Let $\wt\Phi(v):=\Phi(v)-v$. Then, if $\delta$ is sufficiently small, one has
\[
\wt\Phi\ :\ \ol{B}_E(0,2\delta)\to \ol{B}_E(0,\delta),
\]
where $\ol{B}_E$ denotes the closed ball in $E$. Let $f\in B_E(0,\delta)$. Then, the mapping
\[
\wt{\Phi}_f\ :\ v\ \mapsto\ f - \wt\Phi(v)
\]
maps the ball $\ol{B}_E(0,2\delta)$ into itself . Also, $\wt{\Phi}_f$ is continuous in the
weak topology (which is also the weak-$*$ topology, since $E$ is reflexive). Due to the
Banach-Alaouglu theorem (see, e.g.~\cite{RS}~p.~115), $\ol{B}_E(0,2\delta)$ is a compact set
in this topology. Moreover, $\ol{B}_E(0,2\delta)$ is convex and $E$ equipped with the weak
topology is a locally convex space (see, e.g. \cite{RS}~Chapter~V). Therefore, by the
Leray-Schauder-Tyhonoff fixed point theorem (see, e.g.~\cite{RS}~p.~151), there exists $v\in
\ol{B}_E(0,2\delta)$ such that $\wt\Phi_f (v)=v$, i.e., $\Phi (v)=f$.
\end{proof}

\subsection{Local surjection. $\bm{L^p}$ potentials, $\bm{p}\bm{>}\bm{1}$.}
\label{SectLocalSurjLp} Recall that
\[
\cF_\mathrm{cos}: v\mapsto (\nh v^{(0)}; \{-\nh v^{(cn)}\}_{n=1}^{\infty})
\]
is (up~to~a~sign) the cosine-Fourier transform. Let $\cM = \cF_\mathrm{cos} +\wt{\cM}$, where
\[
\wt \cM: v\mapsto (0\,;\{\wt\mu_n(v)\}_{n=1}^{\infty}):=(0\,;\{\mu_n(v)+\nh
v^{(cn)}\}_{n=1}^{\infty})\,.
\]
Let
\[
\cF_\mathrm{cos}^{-1}:(a_0,a_1,\dots)\mapsto a_0 - 2\sum_{n=1}^{\infty} a_n\cos (2\pi nx)
\]
denote the (formal) inverse mapping to $\cF_\mathrm{cos}$.

\begin{proposition}
\label{LocalSurjP>1} Let $1<p<\infty$. Then,

\smallskip

\noindent (i)\phantom{i} the (nonlinear) mapping $\cF_\mathrm{cos}^{-1}\cM$ maps the space
$L^p_\mathrm{even}(0,1)$ into itself;

\smallskip

\noindent (ii) the image $(\cF_\mathrm{cos}^{-1}\cM)(L^p_\mathrm{even}(0,1))$ contains some
neighborhood of~$0$.
\end{proposition}

\begin{proof}
(i) It follows from (\ref{Lasympt}) and $\|v\|_1\le\|v\|_p$ that
\[
\wt \cM \left(L^p_\mathrm{even}(0,1)\right)\ss \wt \cM
\left(L^1_\mathrm{even}(0,1)\right)\ss\ell^{\min\{2,q\}},\qquad q=p/(p\!-\!1),
\]
since  $\wt{\mu}_n(v)=O(n^{-1}\|v\|_1^2\cdot e^{\|v\|_1})$. Thus, the Hausdorff-Young
inequality gives
\[
(\cF^{-1}_\mathrm{cos}\wt \cM)(L^p_\mathrm{even}(0,1))\ss
L^{\max\{2,p\}}_\mathrm{even}(0,1)\ss L^p_\mathrm{even}(0,1)\,.
\]
Moreover, for some constant $C(p)$, one has
\begin{equation}
\label{MapBound} \|(\cF^{-1}_\mathrm{cos}\wt{\cM})(v)\|_p\le C(p)\cdot \|v\|_1^2\cdot
e^{\|v\|_1}.
\end{equation}

\noindent (ii) We are going to apply Lemma~\ref{LocSurjLemma} to the mapping
$\cF_\mathrm{cos}^{-1}\cM$. Note that, for \mbox{$1\!<\!p\!<\!\infty$},
$L^p_\mathrm{even}(0,1)$ is a reflexive Banach space. Due~to~(\ref{MapBound}),
$\cF_\mathrm{cos}^{-1}\cM$ is differentiable (in the Fr\'echet sense) at $0$ and
$d_0(\cF_\mathrm{cos}^{-1}\cM)=I$. This gives the assumption (a1) of Lemma~\ref{LocSurjLemma}.
Thus, it is sufficient to check the assumption (a2), i.e. the continuity of
$\cF_\mathrm{cos}^{-1}\cM$ (or, equivalently, $\cF_\mathrm{cos}^{-1}\wt{\cM}$) in the weak
topology.

Let $v_s\to v$ weakly in $L^p_\mathrm{even}(0,1)$. Let $u\in L^q_\mathrm{even}(0,1)$ and
$h=(h_0,h_1,\dots):=\cF_\mathrm{cos} u$. Then, in order to prove that
$(\cF_\mathrm{cos}^{-1}\wt{\cM})(v_s)\to(\cF_\mathrm{cos}^{-1}\wt{\cM})(v)$ weakly in
$L^p_\mathrm{even}(0,1)$, one needs to show that
\[
\int_0^1
\left((\cF_\mathrm{cos}^{-1}\wt{\cM})(v_s)-(\cF_\mathrm{cos}^{-1}\wt{\cM})(v)\right)\!\!(t)u(t)
dt = 2\sum_{n=1}^{\infty} (\wt{\mu}_n(v_s)\!-\!\wt{\mu}_n(v))h_n\ \to\ 0.
\]
Note that $h\in \ell^{\max\{2,p\}}$ by the Hausdorff-Young inequality, and the norms of the
sequences $\{\wt{\mu}_n(v_s)-\wt{\mu}_n(v)\}_{n=1}^{\infty}$ in $\ell^{\min\{2,q\}}$ are
uniformly bounded due to (\ref{Lasympt}). Thus, Lemma~\ref{WeakCont} and the dominated
convergence theorem imply the result.
\end{proof}

\subsection{Local surjection. $\bm{L^1}$ potentials}
\label{SectL1}

The core argument given in Lemma~\ref{LocSurjLemma} doesn't work for the space $L^1$ since
this space is not reflexive (and is not equipped with any weak-$*$ topology). Nevertheless,
the main result still holds true and the most part of the proof still works well. We start
with some modification of Lemma~\ref{LocSurjLemma}.

\begin{lemma}
\label{LocSurjLemma1} Let $E$ be a Banach space and $\Phi:B_E(0,r)\to E$. Let $F\subset E$ be
a reflexive Banach space and $\|v\|_E\le c\cdot \|v\|_F$ for any $v\in F$ and some constant
$c>0$. If

\noindent \phantom{and} (a1) $\Phi$ is such that $\Phi(v)\!-\!v\in F$ for any $v\in E$ and,
moreover,
\[
\|\Phi(v)-v\|_F=o(\|v\|_E)~\mathit{as}~\|v\|_E\to 0;
\]

(a2) $\Phi$ is continuous in the weak $F$-topology, i.e., for any $v\in E$ and $v_s\!-\!v\in
F$,
\[
v_s\!-\!v\to 0~\mathit{weakly~in~}F~~~\Rightarrow~~~\Phi(v_s)\!-\!\Phi(v)\to
0~\mathit{weakly~in~}F,
\]

\noindent then $\Phi$ is a local surjection, i.e., $\Phi(B_E(0,r))\supset B_E(0,\delta)$ for
some $\delta>0$.
\end{lemma}

\begin{proof} Let $\wt\Phi(v):=\Phi(v)-v$.
It follows from (a1) that, if $\delta$ is sufficiently small,
\[
\wt\Phi\ :\ \ol{B}_E(0,(c\!+\!1)\delta)\to \ol{B}_F(0,\delta).
\]
Let $f\in E$, $\|f\|_E<\delta$, the mapping $\wt\Phi_f:v\mapsto f\!-\!\wt\Phi(v)$ be defined
as in Lemma \ref{LocSurjLemma}, and
\[
\ol{B}_F(f,\delta)\ :=\ \{v\in E: v\!-\!f\in \ol{B}_F(0,\delta)\}.
\]
(note that $f$ and $v$, in general, don't belong to $F$). Since
$\ol{B}_F(0,\delta)\ss\ol{B}_E(0,c\delta)$, one has $\ol{B}_F(f,\delta)\ss
\ol{B}_E(0,(c\!+\!1)\delta)$, and so
\[
\wt\Phi_f\ :\ \ol{B}_F(f,\delta)\to \ol{B}_F(f,\delta).
\]
Moreover, due to (a2), the mapping $\wt\Phi_f$ is continuous in this ``ball'' equipped with
the weak $F$-topology (which is a locally convex topology on this convex compact set). Exactly
as in Lemma~\ref{LocSurjLemma}, the Leray-Schauder-Tyhonoff theorem implies that
$\wt\Phi_f(v)=v$ for some $v\in \ol{B}_F(f,\delta)$, i.e., $\Phi(v)=f$.
\end{proof}

Now we need some modification of Lemma~\ref{WeakCont}, which, together with (\ref{Lasympt}),
implies the assumption (a1) of Lemma \ref{LocSurjLemma1}.

\begin{lemma}
\label{WeakCont1} Let $v_s,v\in L^1(0,1)$ be such that $v_s-v\in L^p(0,1)$ for some
$1<p<\infty$ and $v_s-v\to 0$ weakly in $L^p(0,1)$. Then $\l_n(v_s)\to\l_n(v)$ for any $n\ge
1$.
\end{lemma}

\begin{proof}
Let $u_s:=v_s-v$. Plugging the trivial decomposition $v_s=v+u_s$ into the formula (\ref{VpK=})
for $\vp_k(1,\l,v_s)$, one arrives at
\[
\vp_k(x,\l,q) = \sum_{\{i_1,\dots,i_r\}\ss\{1,\dots,k\}} \iintt_{0\le t_{i_1}\le\dots \le
t_{i_r}\le x} \Phi_\l(t_{i_1},\dots t_{i_r}) \cdot u_s(t_{i_1})\dots u_s(t_{i_r})
dt_{i_1}\dots dt_{i_r}\,,
\]
where
\[
\Phi_\l(t_{i_1},\dots t_{i_r}) = \!\!\!\!\!\!\!\iintt_{0\le t_{j_1}\le\dots \le t_{j_{k-r}}\le
x}\!\! {\textstyle \prod_{m=0}^k \vp_0(t_{m+1}\!-\!t_m,\l)}\cdot v(t_{j_1})\dots
v(t_{j_{k-r}})dt_{j_1}\dots dt_{j_{k-r}},
\]
the sum is taken over all subsets $\{i_1,\dots,i_r\}\ss\{1,\dots,k\}$ of indices and
$\{j_1,\dots,j_{k-r}\}$ denotes the complementary subset. Again, for any fixed $\{i_1,\dots
i_r\}$ and \mbox{$v\in L^1(0,1)$}, the functions
\[
\textstyle f_\l(t_{i_1},\dots,t_{i_r})\ :=\ \chi_{\{0\le t_{i_1}\le \dots \le t_{i_r}\le
x\}}\cdot \Phi_\l(t_{i_1},\dots,t_{i_r}),\qquad |\l|\le M,
\]
form the compact set in $L^q([0,1]^r)$, which gives the result exactly as in
Lemma~\ref{WeakCont}.
\end{proof}

\begin{proposition}
\label{LocalSurjP=1} (i) The mapping $\cF_\mathrm{cos}^{-1}\cM$ maps the space
$L^1_\mathrm{even}(0,1)$ into itself.

\smallskip

\noindent (ii) The image $(\cF_\mathrm{cos}^{-1}\cM) (L^1_\mathrm{even}(0,1))$ contains some
neighborhood of~$0$.
\end{proposition}

\begin{proof}
(i) Recall that (see the proof of Proposition \ref{LocalSurjP>1}(i)) one has
$\cF_\mathrm{cos}^{-1}\cM=I+\cF_\mathrm{cos}^{-1}\wt{\cM}$, and the nonlinear part of our
mapping actually maps $L^1$ potentials into $L^p$ functions (say, for $p=2$,
see~(\ref{MapBound})). In particular, $\cF_\mathrm{cos}^{-1}\cM$ maps the space
$L^1_\mathrm{even}(0,1)$ into itself.

\smallskip

\noindent (ii) Moreover, one has $(\cF_\mathrm{cos}^{-1}\cM])(v)-v\in L^2_\mathrm{even}(0,1)$
for any $v\in L^1_\mathrm{even}(0,1)$ and
\[
\|(\cF_\mathrm{cos}^{-1}\cM)(v)-v\|_2 =o(\|v\|_1)\quad\mathrm{as}\quad\|v\|_1\to 0.
\]
Thus, the assumption~(a1) of Lemma~\ref{LocSurjLemma1} holds with $E=L^1_\mathrm{even}(0,1)$
and $F=L^2_\mathrm{even}(0,1)$. Further, exactly as in Proposition~\ref{LocalSurjP>1}(ii),
Lemma~\ref{WeakCont1} and the dominated convergence theorem give the continuity of the mapping
$\cF_\mathrm{cos}^{-1}\cM$ in the weak $L^2$-topology, i.e., the assumption (a2). So, the
result follows from Lemma~\ref{LocSurjLemma1}.
\end{proof}

\subsection{Global surjection}
\label{SectGlob} To complete the proof of Theorem~\ref{CharEven}, we follow Trubowitz's
approach (cf.~\cite{PT}, pp.~115--116) word for word, if $p>1$, and slightly modify the main
argument, if $p=1$ (cf. the paper~\cite{CKK} devoted to an inverse problem for the perturbed
1D harmonic oscillator, where the same modification was used). The proof of the global
surjection is based on (a)~local surjection near $v=0$ and (b)~explicit solution of the
inverse problem for the perturbation of {\it finitely many} eigenvalues. The latter is given
by
\begin{lemma}[\bf Darboux transform, symmetric case]
\label{DarbouxL} Let $v\in L^1_\mathrm{even}(0,1)$ be a symmetric potential, $n\ge 1$ and $t$
be such that $\l_{n-1}(v)<\l_n(v)+t<\l_{n+1}(v)$. Then there exists a symmetric potential
$v_{n,t}\in L^1_\mathrm{even}(0,1)$ such that
\[
\l_m(v_{n,t})=\l_m(v)\ \mathit{for\ all}\ m\ne n\quad \mathit{and}\quad
\l_n(v_{n,t})=\l_n(v)+t.
\]
Moreover, if $v\in L^p_\mathrm{even}(0,1)$ for some $1\le p<\infty$, then $v_{n,t}\in
L^p_\mathrm{even}(0,1)$ too.
\end{lemma}

\begin{proof}
See \cite{PT}, pp. 107--113, where the modified potential $v_{n,t}$ is constructed explicitly
using the Darboux transform. Namely,
\begin{equation}
\label{V_n_t=}
v_{n,t}=v-2\frac{d^2}{dx^2}\log\{\xi_n(\cdot,\l_n(v)\!+\!t,v);\vp(\cdot,\l_n(v),v)\},
\end{equation}
where $\{f;g\}:=fg'-f'g$ and $\xi_n(\cdot)=\xi_n(\cdot,\l,v)$ denotes the solution of
(\ref{DiffEq}) satisfying the boundary conditions $\xi_n(0)=1$,
$\xi_n(1)=(\vp'(1,\l_n(v),v))^{-1}$ (in particular, the Wronskian is strictly positive on
$[0,1]$). If $v$ is symmetric, then $\xi_n(1)=(-1)^n$ and
$\{\xi_n(\cdot,\l_n(v)+t,v);\vp(\cdot,\l_n(v),v)\}$ is symmetric too. Since
$\{\xi_n;\vp\}'=t\xi_n\vp$, the Wronskian is twice continuously differentiable. In particular,
$v_{n,t}-v_n\in L^p(0,1)$.
\end{proof}

\begin{proof}[\bf Proof of Theorem \ref{CharEven}]
The mapping $\cM$ maps $L^p_\mathrm{even}$ into $\cF_\mathrm{cos}L^p_\mathrm{even}$ (see
Propositions~\ref{LocalSurjP>1}(i),~\ref{LocalSurjP=1}(i)) and is injective due to the well
known uniqueness theorems. Thus, the main problem is to prove that it is surjective. Let
\[
\mu^*=(\mu_0^*,\mu_1^*,\mu_2^*,\dots)\in \cF_\mathrm{cos}L^p_\mathrm{even}(0,1)
\]
be such that $\pi^2+\mu_1^*<4\pi^2+\mu_2^*<\dots$. Since trigonometric polynomials are dense
in $L^p_\mathrm{even}(0,1)$, for any $\delta>0$ there exist some (large) $N$ and a sequence
\[
\mu^\delta=(\mu_0^\delta,\mu_1^\delta,\dots, \mu_N^\delta,\mu_{N+1}^*,\mu_{N+2}^*,\dots)
\]
such that $\pi^2+\mu_1^\delta<4\pi^2+\mu_2^\delta<\dots$ and $\|\cF_\mathrm{cos}^{-1}
\mu^\delta\|_p<\delta$.

Indeed, if $p>1$, then the Fourier series of a function $\cF_\mathrm{cos}^{-1} \mu^*\in
L^p_\mathrm{even}(0,1)$ converge to this function in $L^p$-topology (see,
e.g.~\cite{Edw}~Section~12.10), i.e.,
\[
\|\cF_\mathrm{cos}^{-1} \mu^* - \cF_\mathrm{cos}^{-1}
(\mu_0^*,\mu_1^*,...,\mu_N^*,0,0,...)\|_p=
\|\cF_\mathrm{cos}^{-1}(0,0,...,0,\mu_{N+1}^*,\mu_{N+2}^*,...) \|_p\to 0
\]
as $N\to\infty$, and one can simply take $\mu_0^\delta=\dots=\mu_N^\delta=0$.

If $p=1$, one can still find a finite sequence $(\mu_0^{(N)},\mu_1^{(N)},...,\mu_N^{(N)})$
(or, equivalently, a trigonometric polynomial $2\sum_{n=0}^N\mu_n^{(N)}\!\cos(2\pi n x)$) such
that
\[
\|\cF_\mathrm{cos}^{-1} \mu^* - \cF_\mathrm{cos}^{-1}
(\mu_0^{(N)},\mu_1^{(N)},...\,,\mu_N^{(N)},0\,,0\,,...)\|_1\le\delta,
\]
and take $\mu_n^\delta:=\mu_n^*-\mu_n^{(N)}$ for $j=0,..,N$.

\smallskip

Note that $|\mu_n^\delta|\le\|\cF_\mathrm{cos}^{-1} \mu^\delta\|_1\le\|\cF_\mathrm{cos}^{-1}
\mu^\delta\|_p\le\delta$ for all $n\ge 1$, so the restriction (\ref{MuLe}) holds true. Due to
Proposition~\ref{LocalSurjP>1}(ii) (or Proposition~\ref{LocalSurjP=1}(ii), if $p=1$), there
exists a potential $v^\delta\in L^p_\mathrm{even}(0,1)$ such that $\l_n(v^\delta)=\pi^2n^2+
\mu_0^\delta + \mu_n^\delta$ for all $n\ge 1$, and so
\[
\l_n(v^\delta)=\pi^2n^2+ \mu_0^\delta + \mu_n^*\qquad \mathrm{for\ all}\quad n\ge N\!+\!1.
\]
Adding to $v^\delta$ the constant $\mu_0^*-\mu_0^\delta$ and changing the first $N$
eigenvalues using the procedure given in Lemma \ref{DarbouxL}, one obtains the potential
$v^*\in L^p_\mathrm{even}(0,1)$ such that
\[
\l_n(v^*)=\pi^2n^2+\mu_0^*+\mu_n^*\qquad \mathrm{for\ all}\quad n\ge 1
\]
(to avoid the possible crossing of eigenvalues, i.e., violation of (\ref{MuLe}), during this
procedure, one can always move $\l_1,..,\l_N$ to the far left beginning with $\l_1$ and then
move them to the desired positions beginning with $\l_N$).
\end{proof}

\section{Nonsymmetric case}
\setcounter{equation}{0} \label{SectGen}

\subsection{Preliminaries. Normalizing and norming constants}

If $v$ is not symmetric, then one needs some additional spectral data to determine the
potential uniquely. The possible choices are (cf. \cite{CK09} Appendix B and references
therein):
\begin{itemize}
\item the {\it normalizing constants} (first appeared in Marchenko's paper \cite{Mar})
\[
\a_n(v)= \|\vp(\cdot,\l_n(v),v)\|_2^2=\int_0^1\vp^2(t,\l_n(v),v)dt
=(\dot\vp\vp')(1,\l_n(v),v)\,,
\]
where $\dot\vp$ denotes the derivative with respect to $\l$;
\item the {\it norming constants} introduced by Trubowitz and co-authors (see~\cite{PT})
\[
\n_n(v)=\log[(-1)^n\vp'(1,\l_n(v),v)]\,.
\]
\end{itemize}
Note that
\begin{equation}
\label{a=nu} \a_n(v)=|\dot{w}(\l_n(v),v)| \cdot e^{\n_n(v)},
\end{equation}
where
\[
w(\l,v)\equiv\prod_{m=1}^{\infty}\frac{\l_m(v)-\l}{\pi^2m^2}\,,\qquad \l\in\C,
\]
due to the Hadamard factorization theorem, and so the first factor
\begin{equation}
\label{DotWln=} |\dot{w}(\l_n(v),v)| = \lt|\frac{1}{\pi^2n^2}\prod_{m\ne
n}\frac{\l_m(v)-\l_n(v)}{\pi^2m^2}\rt| = \frac{1}{2\pi^2n^2}\prod_{m\ne
n}\frac{\l_m(v)-\l_n(v)}{\pi^2(m^2-n^2)}
\end{equation}
is uniquely determined by the spectrum. Let
\[
\cF_\mathrm{sin}: v\mapsto \{\nh v^{(sn)}\}_{n=1}^{\infty},\qquad \nh v^{(sn)}=\int_0^1
v(t)\sin (2\pi nt) dt,
\]
be the sine-Fourier transform, and
\[
\cF_\mathrm{sin}^{-1}:(b_1,b_2,\dots)\mapsto 2\sum_{n=1}^{\infty} b_n\sin (2\pi nx)
\]
denote its (formal) inverse. We also use the notation $L^p_\mathrm{odd}(0,1)$ for the space of
all anti-symmetric (or {\it odd}) potentials $v(x)\equiv -v(1\!-\!x)$, $x\in [0,1]$, from
$L^p(0,1)$.

\subsection{Characterization theorem for norming constants}

\begin{theorem}
\label{CharTrubowitz} Let $1\le p<\infty$. The mapping
\begin{equation}
\label{MNmap} \textstyle v\mapsto \left(\cM(v);\cN(v)\right),\qquad \cN(v):=\{2\pi n\cdot
\nu_n(v)\}_{n=1}^{\infty}\,,
\end{equation}
is a bijection between the space of potentials $L^p(0,1)$ and the set of spectral data
$\cM(L^p_\mathrm{even}(0,1))\ts \cF_\mathrm{sin}L^p_\mathrm{odd}(0,1)$. In other words, the
norming constants $\n_n(v)$ multiplied by $n$ can form an arbitrary sequence in
$\cF_\mathrm{sin}L^p_\mathrm{odd}(0,1)$, while the characterization of the possible spectra is
the same as in Theorem~\ref{CharEven}.
\end{theorem}

\begin{proof}
The uniqueness (i.e., the fact that (\ref{MNmap}) is a $1$-to-$1$ map) theorem is well known
(see, e.g.~\cite{PT} p.~62). Further, it directly follows from (\ref{Vp=Series}),
(\ref{VpkEstimate}) that
\[
\n_n(v)=\frac{1}{2\pi n}\cdot \nh v^{(sn)} + O\lt(\frac{\|v\|_1^2e^{\|v\|_1}}{n^2}\rt)\,.
\]
In particular, (\ref{MNmap}) maps $L^p(0,1)$ into $\cM(L^p_\mathrm{even}(0,1))\ts
\cF_\mathrm{sin}L^p_\mathrm{odd}(0,1)$. Moreover, each $\nu_n(v)$ is a continuous function of
the potential in the same sense as in Lemma~\ref{WeakCont}. Repeating the proof of
Proposition~\ref{LocalSurjP>1} (or Proposition~\ref{LocalSurjP=1}, if $p=1$) word for word,
one obtains that (\ref{MNmap}) is a local surjection near $v=0$. Finally (exactly as in
Theorem~\ref{CharEven}), the proof of the global surjection can be finished changing a finite
number of spectral data, which is given by the application of the next (explicit) lemma step
by step.
\end{proof}

\begin{lemma}[\bf Darboux transform, general case]
(i) Let $v\in L^1(0,1)$, $n\ge 1$ and $\l_{n-1}(v)<\l_n(v)+t<\l_{n+1}(v)$. Then there exists a
potential $v_{n,t}\in L^1(0,1)$ such that
\[
\l_m(v_{n,t})=\l_m(v)+t\delta_{nm}\qquad \mathit{and}\qquad \n_m(v_{n,t})=\n_m(v)\
\mathit{for\ all}\ m\ge 1.
\]

\smallskip

\noindent (ii) Let $v\in L^1(0,1)$, $n\ge 1$ and $t\in\R$. Then there exists $v_n^t\in
L^1(0,1)$ such that
\[
\l_m(v_n^t)=\l_m(v)\qquad \mathit{and}\qquad \n_m(v_n^t)=\n_m(v)+t\delta_{nm}\ \ \mathit{for\
all}\ \ m\ge 1.
\]
Moreover, if $v\in L^p(0,1)$ for some $1\le p<\infty$, then $v_{n,t},v_n^t\in L^p(0,1)$ too.
\end{lemma}

\begin{proof}
See \cite{PT}, pp. 91--94, 107--113. The explicit formula for $v_{n,t}$ is given by
(\ref{V_n_t=}) and
\[
v_n^t(x)=v(x)-2\frac{d^2}{dx^2}\log\lt(1-(e^t\!-\!1)\int_x^1\psi_n^2(t,v)dt\rt),
\]
where $\psi_n(\cdot,v)$ is the $n$-th normalized eigenfunction.
\end{proof}

\subsection{Characterization theorem for normalizing constants}

\begin{theorem}
\label{CharMarchenko} Let $1\le p<\infty$. The mapping
\[
\textstyle v\mapsto \left(\cM(v)\,;\,\cA(v)\right),\qquad \cA(v):=\{\pi n\cdot
\log[2\pi^2n^2\a_n(v)]\}_{n=1}^{\infty}\,,
\]
is a bijection between the space of potentials $L^p(0,1)$ and $\cM(L^p_\mathrm{even}(0,1))\ts
\cF_\mathrm{sin}L^p_\mathrm{odd}(0,1)$.
\end{theorem}

\begin{proof} Due to (\ref{a=nu}), (\ref{DotWln=}) and Theorem~\ref{CharTrubowitz}, it is
sufficient to check that
\[
\lt\{ \pi n\cdot \log \prod_{m\ne n}\frac{\l_m(v)\!-\!\l_n(v)}{\pi^2(m^2-n^2)}
\rt\}_{n=1}^{\infty}\in \cF_\mathrm{sin}L^p_\mathrm{odd}(0,1).
\]
Since $\m_m(v)$ are bounded,
\[
\log \frac{\l_m(v)-\l_n(v)}{\pi^2(m^2\!-\!n^2)} = \log \lt( 1+
\frac{\mu_m(v)-\mu_n(v)}{\pi^2(m^2\!-\!n^2)}\rt) = \frac{\mu_m(v)-\mu_n(v)}{\pi^2(m^2-n^2)} +
O\lt(\frac{1}{(m^2\!-\!n^2)^2}\rt).
\]
Summing up over $m\ne n$ (and taking into account that $\mu_n(v)=O(1)$), one obtains
\[
\pi n\cdot \log \prod_{m\ne n}\frac{\l_m(v)\!-\!\l_n(v)}{\pi^2(m^2-n^2)}
=\frac{1}{2\pi}\lt(\sum_{m\ne n}\lt(\frac{1}{m\!-\!n}-\frac{1}{m\!+\!n}\rt)\mu_m(v) -
\frac{1}{2n}\,\mu_n(v) \rt)+ O\lt(\frac{1}{n}\rt)\,.
\]
The error terms belong to $\cF_\mathrm{sin}L^p_\mathrm{odd}(0,1)$ by the Hausdorff-Young
inequality. Denote
\[
\textstyle f:=\cF_\mathrm{cos}^{-1}\left(0,\{\mu_m(v)\}_{m=1}^{\infty}\right)=
-2\sum_{m=1}^{\infty}\mu_m(v)\cos (2\pi mx).
\]
Then, simple straightforward calculations give
\[
\lt\{\frac{1}{2\pi}\lt(\sum_{m\ne n}\lt(\frac{1}{m\!-\!n}-\frac{1}{m\!+\!n}\rt)\mu_m(v) -
\frac{1}{2n}\,\mu_n(v)\rt) \rt\}_{n=1}^{\infty} = \cF_\mathrm{sin}[(\tfrac{1}{2}\!-\!x)f] \,.
\]
Since $(0,\{\mu_m(v)\}_{m=1}^{\infty})\in \cF_\mathrm{cos} L^p_\mathrm{even}(0,1)$, one has
$\cF_\mathrm{sin}[(\tfrac{1}{2}\!-\!x)f]\in \cF_\mathrm{sin}L^p_\mathrm{odd}(0,1)$.
\end{proof}

\end{document}